\documentclass{article}[11pt]
\usepackage{amssymb}
\usepackage{amsmath}
\usepackage{graphicx}
\usepackage{epsfig}
\usepackage[koi8-r,cp1251]{inputenc}
\usepackage[english]{babel}
\title{About statistics of periods of continued fractions of quadratic
irrationalities}
\author{E.\,Yu.~Lerner\footnote{Kazan State University, Russia;
e-mail: eduard.lerner@gmail.com}}
\date{}
\begin{document}
\maketitle
\begin{abstract}
In this paper we answer certain questions posed by V.I.~Arnold,
namely, we study periods of continued fractions for solutions of
quadratic equations in the form $x^2+p x=q$ with
integer~$p$~and~$q$, $p^2+q^2\le R^2$. Our results concern the
average sum of period elements and Gauss--Kuzmin statistics as
$R\to\infty$.
\end{abstract}

\medskip


\section{Statement and discussion of obtained results}
Any value $x\in {\mathbb R}$ is representable as a
continued fraction (CF)
$$
x=a_0+\cfrac{1}{a_1+\cfrac{1}{a_2+\ldots}} =[a_0;a_1,a_2,\ldots],
$$
where $a_0\in {\mathbb Z}$ (if $a_0=0$, then we omit it), $a_i\in
{\mathbb N}$ for all $i\ge 1$. Let $i$ be a fixed position in the
CF expansion and let $P_i(k)$ stand for the probability that
$a_i=k$ with $x$ randomly chosen in the segment~$[0,1)$.
Hereinafter we understand the random choice as a realization of
the uniform distribution. It is well known that according to the
Gauss--Kuzmin theorem,
\begin{equation}
\label{kuz} \lim_{i\to\infty}P_i(k)=\ln(1+\frac{1}{k(k+2)})/\ln 2.
\end{equation}

The roots of the quadratic equation $x^2+p x=q$, where $p,q\in
{\mathbb Z}$, $\Delta=p^2+4q>0$, are expandable in a periodic CF.
Let $x_+(p,q)=\sqrt{\Delta}/2-p/2$. The fractional part of this
value $\{x_+(p,q)\}=x_+(p,q)-\lfloor x_+(p,q) \rfloor$ has no
pre-period in the CF expansion\footnote{see Remark after the proof
of Lemma~3}: $\{x_+(p,q)\}=[[a_1,\ldots,a_{T(p,q)}]]$. Here
$T(p,q)$ is the length of the period of the CF for $x_+(p,q)$;
denote the period itself $(a_1,\ldots,a_{T(p,q)})$ by~$a(p,q)$.

Let us adduce the Arnold conjecture about the Gauss--Kuzmin
statistics for periods of CF of quadratic
irrationalities~(\cite[problem 1993-11B]{problem00},
\cite{ar1,umn,izvestiya}):

\medskip
\noindent{\bf Conjecture~1.\ \ }{\it Choose an integer point
$(p,q)$ in the circle of radius $R$ with the probability
proportional to~$T(p,q)$. Further, let us randomly choose $a_i$
from the corresponding set $a(p,q)$ (V.I.~Arnold performs this
two-step procedure immediately by the random choice from the union
of sets~$a(p,q)$). Then the probability that $a_i=k$ tends to
$\ln(1+\frac{1}{k(k+2)})/\ln 2$ as $R\to\infty$.}
\medskip

Unfortunately, for some reasons, the proof of this conjecture
performed by V.A.~Bykovskii and his followers was not published.
Below we describe the simple proof of the following theorem:

\medskip
\noindent{\bf Theorem~1.\ \ }{\it Let $w\in (0,1)$ be a fixed
value. Let us randomly choose a fixed integer point $(p,q)$ in the
circle of radius $R$. Further, let us choose from the
corresponding set $a(p,q)$ the value $a_i$ with the probability
proportional to $w^i$. Then the probability that $a_i=k$ tends to
$(1-w)\sum_{i=1}^{\infty} P_i(k) w^{i-1}$ as $R\to\infty$. The
limit of this sum with $w \to 1$ is the Gauss--Kuzmin statistics
$\ln(1+\frac{1}{k(k+2)})/\ln 2$.}
\medskip

The latter part of this assertion, evidently, follows from the
Kuzmin theorem~(\ref{kuz}) and the regularity of the Abel
summation method~\cite{hardy}. We did not succeed to prove a
variant of Theorem~1 with $w$ initially equal to~1.

Note that with the random choice of~$(p,q)$ in the circle of radius~$R$
the number $x_+(p,q)$ is not necessarily an irrational real value.
In accordance with Theorem~1 the probability of this event tends to zero as $R\to\infty$.
Let $\Omega_R$ stand for the set of pairs
$(p,q)$ for which this event does not take place with fixed~$R$.

In papers~\cite{umn,izvestiya} V.I.~Arnold experimentally studied
the value $$\widehat T(R)=\sum_{(p,q)\in \Omega_R}
T(p,q)/|\Omega_R|.$$ He established that
\begin{equation}
\label{sim} \widehat T(R)\sim \text{const}\ R\quad\text{ as
$R\to\infty$.}
\end{equation}

Let $T_0(q)=T(0,q)$ be the period of the CF for the square root
of~$q$. Put $\widehat{T}_0(Q)=\sum_{q=1}^Q T_0(q)/Q$. Some experts
in the number theory studied the mean value of the mentioned kind
experimentally; this enabled them to state the following
conjecture: $\widehat{T}_0(Q) \sim \text{const} \sqrt{Q}
\ln^{-\alpha} (Q)$, where $\alpha >0$ (see \cite{Gus12} and
references therein).

Note that one can easily upper estimate
$\widehat{T}_0(Q)$~\cite{av}:
$$\widehat{T}_0(Q) <
\text{const} \sqrt{Q}.
$$
In \cite{GusMS} (see also \cite{Gus12}) E.P.~Golubeva proves that the left-hand side of this
inequality is asymptotically small in comparison with the right-hand one; moreover, if the
extended Riemann conjecture is true,
then the order of their difference is not less than $\ln(Q)^{\ln 2-\varepsilon}$.
The bound
\begin{equation}
\label{down} \widehat{T}_0(Q)> \text{\rm const}\ \sqrt{Q}
\ln^{-\alpha} (Q)
\end{equation}
would allow us to answer the famous Gauss question~\cite{venkov}
about the growth of the number of classes of real quadratic
fields. The known lower bound for $\widehat{T}_0(Q)$ is far from
the right-hand side of~(\ref{down}).
In~\cite{Gus12} E.P.~Golubeva
proves only that $\widehat{T}_0(Q) > \text{const} \ln(Q)$.

Problem~5 listed in~\cite{probl08} implies the estimation of the
growth rate of elements of the period of a CF. Let us adduce its
exact statement~\cite[page~7]{probl08}.

Let $\widehat a (p,q)=\sum_{a_i\in a(p,q)} a_i/T(p,q)$,
$A(R)=\sum_{(p,q)\in \Omega_R} \widehat a (p,q)/|\Omega_R|$. ``The
problem is to evaluate the growth rate of $A(R)$: is $A$ greater
than $C R^\alpha$ for some positive $C$, $\alpha$\,? Or is it
smaller than some $C (\ln R)^\alpha$\,?''.

Let us prove the following assertion.

\medskip
\noindent{\bf Theorem~2.\ \ }{\it
Let us choose $a_i$ in the same way as in the statement of Conjecture~1;
let $A'_R$ stand for the mean value of the chosen number, that is,
\begin{equation}
\label{frac}
A'_R=\frac{\sum_{(p,q)\in\Omega_R}\sum_{a_i\in
a(p,q)}a_i} {\sum_{(p,q)\in\Omega_R} T(p,q)}.
\end{equation}
If correlation (\ref{sim}) or at least that $\widehat T(R)>
\text{\rm const}\ R\, (\ln R)^{-\beta}$ analogous to
inequality~(\ref{down}) takes place, then we have $$A'_R< C (\ln
R)^\alpha.$$}

The main idea of the proof of Theorem~1 is the representation of
the probability under consideration as the Riemann integral. By
applying the Weyl theorem\footnote{The reference to the Weyl
theorem was done by a reviewer of the journal ``Functional
Analysis and Its Application'', the proof of the first variant of
Theorem~1 adduced in~\cite{first} is more awkward.} to the
uniformly distributed sequence $x_+(p,q)$ in Section~2 we easily
prove the desired assertion.

The idea of the application of the Riemann integral sums was used
implicitly for the proof of a similar correlation in the case of
rational values with fixed
denominator~\cite[Theorem~4.5.3.E]{knuth}. However, in this case
one usually applies another technique. An analog of Conjecture~1
for rational values $p/q$ was proved by M.O.~Avdeeva and
V.A.~Bykovskii in~\cite{bik}. The asymptotic validity of the
Gauss--Kuzmin statistics in the case of a fixed denominator
follows from results obtained by H.~Heilbronn and J.W.~Porter (see
\cite{UstinovAlgebraAnaliz}). In ~\cite{us} A.V.~Ustinov
establishes the limit statistics for finite fractions whose
divisors and denominators belong to an arbitrary expanding domain.

The proof of Theorem~2 is based on the upper estimate for
$\sum_{a_i\in a(p,q)} a_i$ as a certain simple function of the discriminant~$\Delta$.
We prove this inequality (possibly, known by experts) with the help of a simple geometric
construction that represents a symbiosis of the ``Nose-Hoover'' algorithm~\cite{ar1} and an explicit technique
for defining values of an integer-valued sign-indefinite quadratic form~\cite{conway}.

Note that sign-indefinite quadratic forms correspond to hyperbolic
operators. Symmetric properties of hyperbolic operators were
studied in geometric terms by F.~Aicardi (\cite{AicardiFAOM}, see
also \cite{AicardiArchiV,Karpenkov}). The algebraic approach to
this topic is described in well-known
papers~\cite{Zagier,ManinMarco}.

This paper is written under the bright impression of lectures
delivered by V.I.~Arnold~\cite{vid1,vid2,vid3} and a result of the
performed numerical experiments.

\section{Proof of Theorem~1.}

Let $b$ be an arbitrary number from the segment $[0,1]$.
Recall~(\cite{keipers}) that a sequence is said to be uniformly
distributed, if among its first $N$ elements the quotient of terms
whose fractional part is less than $b$ tends to $b$ as
$N\to\infty$.

\medskip

\noindent{\bf Lemma 1.}\ \ {\it Let us enumerate numbers $x_+(p,q)$,
$(p,q)\in \mathbb{Z}^2$, $p^2+4q>0$ in ascending order of the
distance from the point $(p,q)$ to the origin of coordinates
(in the case of equal distances the numeration is arbitrary). The obtained
sequence is uniformly distributed.}

\medskip

Since $x_+(p,q)-x_+(-p,q)\in\mathbb{Z}$, in order to prove the lemma, suffice
it to consider the terms of the sequence
that correspond to points in the right half-plane. Lemma~1 easily follows from Lemma~2.

\medskip

\noindent{\bf Lemma 2.}\ \ {\it Let $n$ is a fixed integer number.
The subsequence of the sequence from Lemma~1 that corresponds to
points $(p,q)$ such that $n\le q/p<n+1$, $p>0$ is uniformly
distributed.}

\medskip

\begin{figure}
\begin{center}
\includegraphics{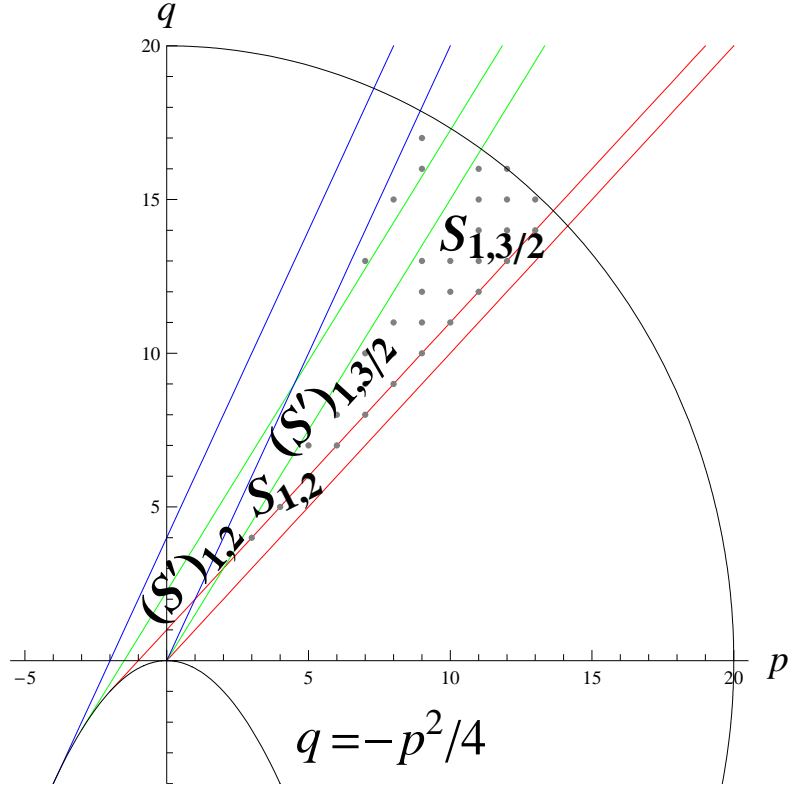}
\end{center}
\caption{Domains considered in the proof of Lemma~2.
The rays that bound the domain $S'$ are tangent lines to the parabolic curve
$q=-p^2/4$. Indicated are sets of points $M_{20,1,1/2}\cap
M'_{20,1,1/2}$ and $\bar M_{20,1,1/2}\cap \bar M'_{20,1,1/2}$.}
\label{Pic1}
\end{figure}

\noindent{\it Proof of Lemma 2}\ \ Conditions $n\le q/p<n+b$,
$p>0$, $b\in [0,1]$ define the sector $S_{n,n+b}$ in the plane
$(p,q)$. Instead of it, let us consider another
domain~$S'_{n,n+b}$. The domain~$S'_{n,n+b}$ is defined by
inequalities $n p+n^2\le q< (n+b) p+(n+b)^2$, with $p>-2n$ and
inequalities $-p^2/4\le q<(n+b)p+(n+b)^2$ with $-2n-b<p\le -2n$
(see Fig.~\ref{Pic1}). One can easily make sure that the rays that
define the boundaries of the domain~$S'$ are equiscalar lines of
the function~$x_+(p,q)$. Evidently, the condition $\{x_+(p,q)\}<b$
for points from $S'_{n,n+1}$ is equivalent to their belonging to
the subsector~$S'_{n,n+b}$.

Let $M_{R,n,b}$ ($M'_{R,n,b}$) stand for the set of integer points
located inside the sector $S_{n,n+b}$ (the domain $S'_{n,n+b}$)
and the circle of radius~$R$ centered at the origin of
coordinates. Let $\bar M_{R,n,b}=M_{R,n,1}\setminus M_{R,n,b}$,
$\bar M'_{R,n,b}=M'_{R,n,1}\setminus M'_{R,n,b}$.

Evidently, for $R\to\infty$ we have $|M'_{R,n,b}\cap
M_{R,n,b}|/|M_{R,n,1}|\to b$, $|\bar M'_{R,n,b}\cap \bar
M_{R,n,b}|/|M_{R,n,1}|\to 1-b$, whence we obtain the assertion of
Lemma~2. \qquad $\Box$

The assertion of Lemma~1 follows from the inequality $|\bigcup_{n:|n|\le N}
M_{R,n,1}|/|M_{R,+}|>1-1/N$, where $M_{R,+}$ is the collection of all integer points
located in the half-circle of radius~$R$.

Let us now immediately prove Theorem~1. Let $i,k\in \mathbb{N}$,
let $I^i_k(x)$ stand for the indicator function of the set of
points~$x$ from the segment~$[0,1]$ for which the $i$-th position
in the CF is occupied by the number~$k$. This set and its
complement are representable as unions of countable numbers of
intervals~\cite{hinchin}. Consequently, the function $I^i_k(x)$ is
countably continuous and therefore it is Riemann integrable. The
function $f_{k,w}(x)=(1-w)\sum_{i=1}^\infty I^i_k(x) w^{i-1}$ has
the same property. The integral of this function
\begin{equation}
\label{int} \int_{0}^1 f_{k,w}(x)\,dx=(1-w)\sum_{i=1}^{\infty}
P_i(k) w^{i-1}
\end{equation}
equals the probability of the event~$A$: $a_j=k$ with the random
choice of the number $x$ in the segment $[0,1)$. Here $j$ is
chosen randomly, namely, $\text{Prob}(j=i)=(1-w) w^{i-1}$,
$i=1,2,3,\ldots$.

We can calculate the Riemann integral~(\ref{int}) with the help of
the Weyl theorem~\cite{keipers}. In accordance with this theorem
the frequency of the event~$A$ for the first $N$ terms of a
uniformly distributed sequence tends to its probability. Taking
into account the fact that the CF for $\{x_+(p,q)\}$ has no
pre-period, we obtain the assertion of Theorem~1.\qquad $\Box$

\section{Proof of Theorem~2}

We remind~\cite{ar1} that the ``Nose-Hoover'' algorithm for
finding a CF of a real number $x$, $x>0$, is reduced to the
geometric method which constructs the boundary of the convex shell
of the set of integer nonnegative points~$(u,v)$ located above
(below) the straight line~$v=x u$. Let $\pmb{e}_0$ stand for the
vector $(0,1)$, $\pmb{e}_1=(1,0)$. Put
\begin{equation}
\label{step} \pmb{e}_{n+1}=\pmb{e}_{n-1}+a_{n-1} \pmb{e}_n,\quad
n=1,2,\ldots,
\end{equation}
where $a_{n-1}$ is the maximal integer such that the vector
$\pmb{e}_{n+1}$ lies below (for odd $n$) or above (for even $n$)
the straight line~$v=x u$. Klein noted that $a_n$, $n=0,1,\ldots$,
coincide with partial quotients of the CF of the number~$x$. The
geometric algorithm results in two polylines which represent parts
of sails of the CF.

Let $x>1$. We need a slightly modified algorithm which results in
one infinite polyline~$L$, originating at the point~$(1,1)$, whose
segments are the vectors $a_{n-1} \pmb{e}_{n}$, $n=1,2,\ldots$.
Evidently, $\pmb{e}_{n+1}+\pmb{e}_n=\pmb{e}_n+\pmb{e}_{n+1}$.
Consequently, if in the standard ``Nose-Hoover'' algorithm
~(\ref{step}) we add the extra vector $\pmb{e}_n$ and thus go out
of the line $v=x u$ (i.\,e., we consider the vector
$\pmb{e}_{n-1}+(a_{n-1}+1) \pmb{e}_n$), then we get the first
integer point on the segment of the polyline of another sail
constructed at step $n+1$. This fact justifies a simple geometric
algorithm which constructs the polyline~$L$.

\begin{figure}
\begin{center}
\includegraphics{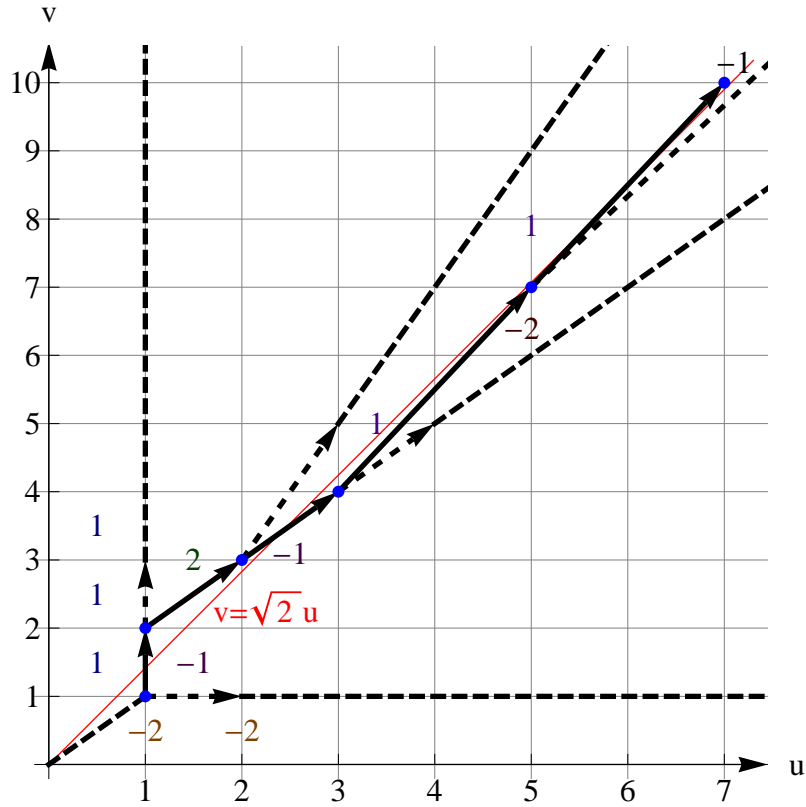}
\end{center}
\caption{The gradual ``Nose-Hoover'' algorithm for $x=\sqrt{2}$}
\label{Pic2}
\end{figure}

We begin the construction process with the point $(1,1)$; it is
convenient to connect it with the origin of coordinates by the
segment which does not enter in~$L$. As the ``constructive'' term
at the first ($n$th) step we choose the vector $\pmb{e}_1$ (the
vector $\pmb{e}_n$). We add this vector till we go out of the line
$v=x u$. The newly added vector $\pmb{e}_{n+1}$ is directed from
the origin of coordinates to the point obtained as a result of the
latter addition up to the step out of the line. See
Fig.~\ref{Pic2} for the first segments of the polyline for
$x=\sqrt{2}$. The zero partial quotient (i.\,e., $\lfloor x
\rfloor$) differs from zero. In order to provide the
correspondence to the indices of partial quotients, it is also
convenient to begin the numeration of segments of the polyline~$L$
with zero.

Let us consider the approximation of a number~$x$ by mediants. Let
$f_1=\frac{v_1}{u_1}$, $f_2=\frac{v_2}{u_2}$ be irreducible
fractions such that $v_1,u_1,v_2,u_2\ge 0$ and $f_1<f_2$. Let us
define the operation $\downarrow$ of finding the mediant (the
``insertion'' operation) by the formula $f_1\downarrow
f_2=\frac{v}{u}$, where $v=v_1+v_2$, $u=u_1+u_2$. Earlier we used
the denotation introduced by A.A.~Kirillov. Evidently, in the
geometric representation of the fraction $\frac{v}u$ as the
integer vector $(u,v)$ the operation $\downarrow$ corresponds to
the addition of vectors $(u_1,v_1)$ and $(u_2,v_2)$. The obtained
diagonal of the parallelogram appears to be ``inserted'' between
its sides.

Let $x\in (f_1,f_2)$, $f_3=f_1\downarrow f_2$. Since $f_3\in
(f_1,f_2)$, one can treat $f_3$ as an approximation of the
number~$x$. If $x\in (f_1,f_3)$, then we put $f_4=f_1\downarrow
f_3$, otherwise we do $f_4=f_3\downarrow f_2$. The process of the
approximation of the number~$x$ by mediants consists in the
repetition of these operations for the corresponding intervals.

The condition $x>0$ means that $x\in\left( \frac01,\frac10
\right)$. The initial approximation of $x$ for such an interval is
the fraction~$\frac11$. One can easily see that for the mapping
$$ \text{the fraction $\frac{v}u \quad \leftrightarrow\quad$ the point
$(u,v)$}
$$
the geometric representation of the algorithm for approximating
the number~$x$ by mediants is completely identical to the
algorithm for constructing the polyline~$L$. The technical
distinction consists in the following fact: earlier we constructed
each segment of the polyline~$L$ ``at once'', but now it ``grows
gradually'' due to the stepwise addition of the next
vector~$\pmb{e}_i$.

Further we consider a simplified version of this algorithm
for the quadratic irrationality $x_+(p,q)$.
We will need it for the proof of the next assertion.

\medskip

\noindent{\bf Lemma~3.}\ \ {\it For $n\in \mathbb{N}$ put
$D(n)=\sum_{u=1}^{\lfloor\sqrt{n}\rfloor}\tau(n-u^2)$, where
$\lfloor\cdot\rfloor$ is the integer part, $\tau(m)$ is the
quantity of divisors of the number~$m$. We have the inequality
\begin{eqnarray}
\label{th2} \sum_{a_i\in a(p,q)} a_i\le f(\Delta/4),\qquad
\text{ where\quad $\Delta=p^2+4q$,\phantom{nbbnbnbnb}}\\
\label{th2even} f(n)= 2 D(n)+\tau(n), \text{if $n$ is integer;}\\
\label{th2odd}
f(n)=2\mathop{\sum\limits_{i\in\{1,3,\ldots\}}}\limits_{i^2<4n}\tau(n-i^2/4),
\text{  otherwise.   }
\end{eqnarray}
In the case of an odd period~$T(p,q)$ one can improve this bound,
namely, divide the right-hand side of inequality~(\ref{th2}) by
two.}
\medskip

\noindent{\it Proof of Lemma~3}\ \ Note that
$\{x_+(p,q)\}=\{x_+(p+2,q-p-1)\}$, therefore, without loss of
generality, we assume that $p=0$ or $p=1$, here $x_+(p,q)>1$ (the
special case $x_+(1,1)$ is evident) and the second root of the
quadratic equation is negative.

Let us apply to $x_+(p,q)$ the gradual ``Nose-Hoover'' algorithm
described above. Let $(u_n,v_n)$ stand for the $n$th integer point
on the polyline~$L$, $n=0,1,2,\ldots$. We obtain it at the $n$th
step of the gradual ``Nose-Hoover'' algorithm. In Fig.~\ref{Pic2},
for example, $(u_0,v_0)=(1,1)$, $(u_1,v_1)=(1,2)$,
$(u_2,v_2)=(2,3)$, $(u_3,v_3)=(3,4)$, $(u_4,v_4)=(5,7)$, etc.

Let us draw one more vector which originates from the point
$(u_n,v_n)$; denote it by $\pmb{e}'_n$. The vector $\pmb{e}'_n$ is
defined by the condition $\pmb{e}'_n+\pmb{e}_i+\pmb{e}''_n=0$,
where $\pmb{e}_i$ is the vector which originates from the point
$(u_n,v_n)$ and goes along the polyline~$L$,
$\pmb{e}''_n=-(u_n,v_n)$ (we considered this vector earlier). The
collection of vectors $\{\pmb{e}'_n,\pmb{e}_i,\pmb{e}''_n\}$,
where each one is defined accurate to the multiplier $(-1)$, is a
{\it superbasis}~\cite{conway}. This means that any pair of these
vectors generates the whole integer lattice. The transition from
the point $(u_n,v_n)$ to that $(u_{n+1},v_{n+1})$ corresponds to
the replacement of one superbasis with another one; the latter
differs from the initial superbasis only in one of three its
elements. Thus, for the polyline~$L$ in Fig.~\ref{Pic2} the
initial superbasis is $\{(0,1),(1,0),(1,1)\}$, the 1st one is
$\{(1,1),(1,0),(1,2)\}$, the 2nd one is $\{(1,2),(1,1),(2,3)\}$,
the 3rd one is $\{(1,1),(2,3),(3,4)\}$, etc.

Let us extend vectors~$\pmb{e}'_n$ up to the rays which originate
at the points $(u_n,v_n)$. As a result, the first quadrant appears
to be divided onto disjoint connected domains (see
Fig.~\ref{Pic2}). Let us associate each vector $(u_n,v_n)$ with a
domain, whose boundary contains all points of polylines which
include the vector $(u_n,v_n)$ in the corresponding superbasis. In
the picture these domains look like narrow ``crevices'' located at
the north-east of the points $(u_n,v_n)$. We associate the domains
which border on the coordinate axes with the unit vectors of these
axes.

We have
\begin{equation}
\label{equiv} \text{$(u_n,v_n)$ is located above the line
$v=x_+(p,q)u$}\ \Longleftrightarrow\ v_n^2+p v_n u_n>q u_n^2.
\end{equation}

Let us write the values of the quadratic form $v^2+p v u-q u^2$ at
the points $(u_n,v_n)$ in the corresponding domain. See
Fig.~\ref{Pic2} for the result obtained in the case $p=0$, $q=2$
(the upper right corner of the figure contains the value
calculated at the point~$(5,7)$). From~(\ref{equiv}) it follows
that the points of the polyline~$L$ belong to boundaries of
domains with positive and negative values. The collection of
superbases which correspond to these points is called the {\it
river} of the quadratic form~\cite{conway}.

\begin{figure}[t]
\begin{center}
\begin{picture}(100,40)
\put(30,20){\line(-1,1){20}} \put(30,20){\line(-1,-1){20}}
\put(30,20){\vector(1,0){20}} \put(50,20){\line(1,0){20}}
\put(70,20){\line(1,1){20}} \put(70,20){\line(1,-1){20}}
\put(14,18){d} \put(80,16){c} \put(48,37){a} \put(48,3){b}
\put(52,22){h}
\end{picture}
\end{center}
\caption{The arithmetic progression rule} \label{Pic3}
\end{figure}
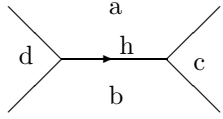

Let us adduce the main result of the theory of quadratic forms
which enables one to make the calculation of the period of a CF
much easier. Using the values of the quadratic form at three
vectors from any superbasis, one can easily restore the CF. To put
it more precisely, the following correlation~\cite{conway} is
true. Let domains associated with values $a,b,c,d$ of the
quadratic form be located in accordance with Fig.~\ref{Pic3}.
Then values $d,a+b,c$ form an arithmetic progression. Moreover, if
$h$ is the common difference of this progression, then
\begin{equation}
\label{main} h^2/4-a b= \Delta/4.
\end{equation}

Therefore, using the values $(a,b)$ for two neighboring domains
with different signs separated by the polyline~$L$, and the value
$h$, one can unambiguously restore all subsequent and previous
values of partial quotients of the CF, moving along the river of
the quadratic form (see the figures in Chapter~1 of
book~\cite{conway}).

\begin{table}[t]
\begin{center}
\begin{tabular}{|c|c|c|c|c|c|}
\hline
{\it step number} &0&1&2&3&4\\
\hline
$a$ is the value located at the north west of~$L$&1&2& 1&1&1\\
\hline
$b$ is the value located at the south east of~$L$&-1&-1&-1&-2&-1\\
\hline
$h$ is the common difference of the arithmetic progression&2&0&-2&0&2\\
\hline
{\it the number of a segment of the polyline } &\it0&\it1&\it1&\it2&\it2\\
\hline
\end{tabular}
\end{center}
\caption{All values of $(a,b,h)$ for the polyline
represented in Fig.~\ref{Pic2}}
\end{table}

Let us now immediately prove the inequalities~(\ref{th2}). They
follow from correlation~(\ref{main}). Assume that at step~$n_1$ we
get the same parameters of the arithmetic progression as those
obtained at step~$0$. Such a step exists because the algorithm is
invertible. Let~$l_1$ stand for the number of the polyline segment
which contains the point $(u_{n_1},v_{n_1})$. Then a part of the
CF $a_{1},a_{2},\ldots,a_{l_1}$ becomes periodic. In addition, the
number $l_1$ is even, because the points $(u_{0},v_{0})$ and
$(u_{n_1},v_{n_1})$ are located to one side of the polyline~$L$.
In the case of an odd period the sequence
$a_{1},a_{2},\ldots,a_{l_1}$ contains at least two repeating
subsequences. In accordance with the algorithm we have
$\sum_{i=1}^{l_1} a_{i}=n_1$. Therefore, suffice it to estimate
the number of all possible triplets $(a,b,h)$ which
satisfy~(\ref{main}). The obtained bound obeys
formulas~(\ref{th2even},\ref{th2odd}). The multiplier~2 in the
right-hand sides of these formulas appears, because we have to
take into account various signs of~$h$, and the term
$\tau(\Delta/4)$ does when we consider the case $h=0$. Lemma~3 is
proved. \qquad $\Box$

\medskip

\noindent{\it Remark}\ \ \ In the proof of Lemma~3 we established
that CF for $\{ x_+(p,q)\}$ has no pre-period.

\medskip

Let us complete the proof of Theorem~2. Formally speaking, the
Dirichlet theorem (the equality $\sum_{u=1}^n \tau(u)/n=\ln
n+(2\gamma-1)+O(1/\sqrt{n})$) does not imply that $D(n)\sim
\sqrt{n} \ln n$, because the values of the function $\tau$ are not
uniform. A more accurate bound is
$$
D(n)=O(\ln^3(n)\sqrt{n})
$$
(see \cite{av} and references therein). This, evidently, implies
that $f(n)=O(\ln^3(n)\sqrt{n})$. Taking the sum over all possible
values of~$(p,q)$, we obtain that the numerator of
fraction~(\ref{frac}) equals $O(R^3 \ln^3(R))$, whereas the
denominator by assumption is greater than $\text{\rm const}\ R^3\,
(\ln R)^{-\beta}$. Theorem~2 is proved.\qquad $\Box$

\end{document}